\documentclass[12pt]{amsart}

\usepackage{amsmath,amssymb,amsthm}
\usepackage[all]{xy}
\usepackage{amscd}
\usepackage{amsfonts}
\usepackage{color}
\usepackage[dvipdfm,colorlinks=true]{hyperref}

\numberwithin{equation}{section}

\SelectTips{eu}{12}

\theoremstyle{definition}

\theoremstyle{remark}

\newcommand{\A}{\underline{A}}
\newcommand{\X}{\underline{X}}
 \newcommand{\fk}{\mathbf{k}}
\newcommand{\link}{\mathrm{link}}
\newcommand{\Sd}{\mathrm{Sd}}
\newcommand{\st}{\mathrm{star}}

\newcommand{\Tor}{\mathrm{Tor}}
\newcommand{\Z}{\mathbb{Z}}

\makeatletter
  \def\wideubar{\underaccent{{\cc@style\underline{\mskip10mu}}}}
  \makeatother

\title{A Golod complex  with non-suspension moment-angle complex}

\author{Kouyemon Iriye}
\address{Department of Mathematics and Information Sciences, Osaka Prefecture University, Sakai, 599-8531, Japan}
\email{kiriye@mi.s.osakafu-u.ac.jp}
\author{Tatsuya Yano}
\address{Department of Mathematics and Information Sciences, Osaka Prefecture University, Sakai, 599-8531, Japan}
\email{su301034@osakafu-u.net}

\subjclass[2010]{55P15(13F55)}
\keywords{Stanley-Reisner ring, Golod, Massey product, Hopf map}

\begin{document}

\maketitle

\begin{abstract}
It could be expected that the moment-angle complex associated with a Golod simplicial complex is homotopy 
equivalent to a suspension space. In this paper, we provide a counter example to this expectation. 
We have discovered this complex through the studies of 
the Golod property of the Alexander dual of a join of simplicial complexes, and
that of a union of simplicial complexes.
\end{abstract}

\baselineskip 16pt

\section{Introduction}

The {\it Stanley-Reisner ring} (or {\it face ring}) of a simplicial complex $K$ over an index set
$[m]=\{1,\cdots,m\}$ is defined as the quotient graded algebra
\[
\fk[K]=\fk[v_1,\cdots,v_m]/\mathcal{I}_K,
\]
where $\fk$ is a commutative ring with unit and
$\mathcal{I}_K=(v_{i_1}\cdots v_{i_k}\ |\ \{i_1,\cdots,i_k\}\not\in K)$ is the
 {\it Stanley-Reisner ideal} of $K$.
$K$ is called {\it Golod} over a field $\fk$ if its Stanley-Reisner ring $\fk[K]$
is Golod over $\fk$. That is, the multiplication and all higher Massey products in
\[
\Tor^{\fk[v_1,\cdots,v_m]}(\fk[K],\fk)=H(\Lambda[u_1,\cdots,u_m]\otimes\fk[K],d)
\]
are trivial, where the Koszul differential algebra
$(\Lambda[u_1,\cdots,u_m]\otimes\fk[K],d)$ is the bigraded differential algebra with
$\deg u_i=(1,2)$, $\deg v_i=(0,2)$, and $du_i=v_i$ for $i=1,\cdots,m$. Originally, the algebra $\fk[K]$ or the
ideal $\mathcal{I}_K$ was defined to be Golod if the following equation holds:
\[
\sum_{i\geq0;\ j\geq0}\dim_{\fk}\Tor^{\fk[K]}_{j,2i}(\fk,\fk)t^jz^i=
\frac{(1+tz)^m}{1-t\sum_{i\geq0;\ j\geq1}\dim_{\fk}\Tor_{j,2i}^{\fk[v_1,\cdots,v_m]}(\fk[K],\fk)t^jz^i},
\]
where $\Tor^{\fk[K]}_{j,2i}(\fk,\fk)$ and $\Tor_{j,2i}^{\fk[v_1,\cdots,v_m]}(\fk[K],\fk)$
denote the homogeneous components of degree $2i$. Golod \cite{G}
proved the equivalence of the two conditions, and thereafter his name has been used to
refer a ring that satisfies the condition. The reader may also refer to Gulliksen and Levin \cite{GL} or  Avramov \cite{Av}.

Baskakov, Buchstaber, and Panov \cite{BBP} and Franz \cite{Fr} independently demonstrated
that the torsion algebra
$\Tor^{\fk[v_1,\cdots,v_m]}(\fk[K],\fk)$ is isomorphic to the cohomology ring of the moment-angle complex $Z_K$
associated with $K$.

\bigskip

{\bf Theorem 1.1(\cite{BBP, Fr}).} {\it For a commutative ring $\fk$ with unit, the following isomorphisms of algebras hold:
\[
H^*(Z_K;\fk)\cong \Tor_*^{\fk[v_1,\cdots,v_m]}(\fk[K],\fk)\cong \bigoplus_{I\subset[m]}\tilde{H}^*(K_I;\fk),
\]
where $\tilde{H}^*(K_I;\fk)$ denotes the reduced cohomology of the full subcomplex $K_I$ of $K$ on $I$, and
$\tilde{H}^*(K_\emptyset;\fk)=0$ for $*\neq-1$ and $=\fk$ for $*=-1$. The last isomorphism is the sum
of isomorphisms given by
\[
H^p(Z_K;\fk)\cong \bigoplus_{I\subset[m]}\tilde{H}^{p-|I|-1}(K_I;\fk),
\]
and the ring structure is given by the maps
\[
\tilde{H}^{p-|I|-1}(K_I;\fk)\otimes \tilde{H}^{q-|J|-1}(K_J;\fk)\to \tilde{H}^{p+q-|I|-|J|-1}(K_{I\cup J};\fk)
\]
that are induced by the canonical inclusion maps $\iota_{I,J}:K_{I\cup J}\to K_I*K_J$ for $I\cap J=\emptyset$ and zero otherwise, where $K_I*K_J$ denotes the join of two simplicial complexes $K_I$ and $K_J$.
}

\bigskip

Here, we recall that if the moment-angle complex $Z_K$ is homotopy equivalent to a suspension, then
the multiplication and all higher Massey products in $H^*(Z_K;\fk)$ are trivial. For example, see
Corollary 3.11 of \cite{P}. That is, the following implication holds:
\begin{equation}
\text{$Z_K$ is homotopy equivalent to a suspension}\Longrightarrow \text{ $K$ is Golod},
\end{equation}
where $K$ is Golod if $K$ is Golod over any field $\fk$.
 This observation enables us to investigate the Golod property through the
study of moment-angle complexes. One of the first studies in this direction was introduced by Grbi\'{c} and
Theriault \cite{GT}. They demonstrated that the moment-angle complex associated with a shifted
simplicial complex is homotopy equivalent to a wedge of spheres.
 In \cite{IK14}, Kishimoto and the first author extended this result to dual sequentially
Cohen-Macaulay complexes, and provided some new Golod complexes. In these studies, the following theorem
concerning the decomposition of polyhedral products (see Definition 2.1),
as introduced by Bahri, Bendersky, Cohen, and Gitler \cite{BBCG}, plays an essential role.

\bigskip

{\bf Theorem 1.2(\cite{BBCG}).} {\it Let $K$ be a simplicial complex on $[m]$ and let 
$(C\X,\X)=\{(CX_i,X_i)\}_{i\in[m]}$, where each $X_i$ is a based space and $CX$ is
the reduced cone of a based space $X$. Then, the following homotopy equivalence holds:
\[
\Sigma Z_K(C\X,\X)\simeq \Sigma\bigvee_{I\subset[m]}\Sigma|K_I|\wedge\widehat{X}^{I},
\]
where $\widehat{X}^{I}=\wedge_{i\in I}X_i$ and $\widehat{X}^\emptyset=*$.
}

\bigskip

We call this decomposition of polyhedral products the BBCG decomposition for $K$. If this decomposition is desuspendable, i.e., if the homotopy equivalence
\[
Z_K(C\X,\X)\simeq \bigvee_{I\subset[m]}\Sigma|K_I|\wedge\widehat{X}^{I}
\]
holds for any sequence of based CW-complexes $\X$, then we say that the BBCG decomposition is 
{\it desuspendable} for $K$. In particular, this implies that $Z_K$ is homotopy equivalent to a suspension. 

\bigskip

In this paper, we study the Golod properties of the Alexander dual of $K*L$ and 
$K\cup_{\alpha}L$, where $\alpha$ is a common face of $K$ and $L$.
The precise statements of the results are given in the next section.

By Theorem 1.1, the multiplicative structure of $H^*(Z_K;\fk)$ is trivial if and only if the maps
$\iota_{I,J}:K_{I\cup J}\to K_I*K_J$ for $I\cap J=\emptyset$ induce the trivial maps on the reduced 
cohomology theory.
By strengthening this condition, $K$ is said to be {\it (stably) homotopy Golod} \cite{IK14} if the maps
$|\iota_{I,J}|:|K_{I\cup J}|\to |K_I*K_J|$ for $I\cap J=\emptyset$ are (stably) null homotopic and
$H^*(Z_K;\fk)$ has trivial higher Massey products for any fields $\fk$.
By definition, the following implication holds:
\[
\text{stably homotopy Golod}\Longrightarrow \text{Golod}.
\]
The second purpose of this paper is to prove that this implication is strict.

\bigskip

{\bf Theorem 1.3.} {\it There is a Golod simplicial complex $K$ such that $K$ is not stably homotopy
Golod. Moreover, $Z_K$ can be chosen to be torsion free.}

\bigskip

Here, a space or a simplicial complex $X$ is called {\it torsion free} if its integral homology groups $H_*(X;\Z)$
are torsion free. In general, $Z_K$ can have torsion in its  integral homology groups even if $K$ is Golod. 
The 6 vertex triangulation of $\mathbb{R} P^2$ is such an example, see for example \cite{GPTW} or 
Example 10.10 of \cite{IK14}.

It could be expected that the converse of the implication (1.1) is also true.
Theorem 1.3 provides a counter example to this expectation.
In fact, if $Z_K$ is homotopy equivalent to a suspension then the fat wedge filtration of $Z_K$ is trivial, 
by Theorem 1.3 of \cite{IK14}. By Theorem 6.9 of the same paper, we see that $K$ is stably homotopy 
Golod, which contradicts our result. Thus, $Z_K$ is not homotopy equivalent to a suspension.

In the next section, we state the main results of this paper. The subsequent sections are devoted to
their proofs.

\bigskip

{\bf Acknowledgments}  We are grateful to Daisuke Kishimoto for providing useful discussions and 
pointing out the paper \cite{MS},
which was applied in the construction of $K$ in Theorem 1.3. We also thank Lukas Katth\"{a}n for
his comments on the first draft of this paper. He kindly pointed out the papers \cite{HS} and \cite{D},
and provided us with ideas for simplifying the arguments in the subsection 5.3. We would also like to thank the 
anonymous referee for his/her careful reading of our manuscript and his/her many constructive comments 
and suggestions, which helped us to improve the manuscript. Finally, the first author thanks Japan Society 
for the Promotion of Science (KAKENHI grant number 26400094) for funding.

\section{Results}

In this section, we state our main results. We begin by setting notation regarding simplicial complexes.

Let $K$ be a simplicial complex on an index set $V$. In this paper, we only consider finite
simplicial complexes. A {\it facet} of a simplicial complex is a maximal face with respect to the inclusion, 
and a subset $\sigma\subset V$ is called a {\it minimal non-face} if $\sigma\not\in K$ 
but $\sigma-v\in K$ for every $v\in\sigma$. The subset $V(K)$ of $V$ defined by
$V(K)=\cup_{\sigma\in K}\sigma$ is called the vertex set of $K$, and an element of $V-V(K)$
is called a ghost vertex. In this paper we allow a simplicial complex to have a ghost vertex. 
But we assume that a simplicial complex has the simplex $\emptyset$. 
 
For a finite set $V$, we denote the full simplex on $V$ by $\Delta^V$.
Its boundary is denoted by $\partial\Delta^V$.
We also use the symbol $\Delta^n$ to denote an $n$-dimensional simplex.
$|K|$ denotes a geometric realization of $K$.
The link and star of a face $\sigma$ of $K$ are denoted by $\link_K(\sigma)$ and
$\st_K(\sigma)$, respectively:
\[
\link_K(\sigma)=\{\tau\subset \sigma^c=V-\sigma\ |\ \tau\cup \sigma\in K\},\ 
\st_K(\sigma)=\{\tau\subset V\ |\ \tau\cup \sigma\in K\}. 
\]
For a subset $I\subset V$, $K_I$ denotes the full subcomplex of $K$ indexed by $I$; that
is, \\
$K_I=\{\sigma\subset I\ |\ \sigma\in K\}$. Then, $K_I$ is called the restriction of $K$ to $I$. 
For a vertex $v$ of K, the deletion of $v$ from K is denoted by $K-v =K_{V-{v}}$. 
For simplicial complexes $K$ and $L$ with disjoint index sets $V$ and $W$,
the simplicial {\it join} $K*L$ on the index set $V\sqcup W$ is defined
by $K*L=\{\sigma\sqcup\tau\subset V\sqcup W\ |\ \sigma\in K,\ \tau\in L\}$,
where $\sqcup$ always denotes the disjoint union of sets. We write $K*c$ to stand for $K*\Delta^{\{c\}}$. 
For disjoint subsets $I,J\subset V$ 
the canonical inclusion map $\iota_{I,J}:K_{I\sqcup J}\to K_I*K_J$ mentioned in Theorem 1.1 is defined 
as $\iota_{I,J}(\sigma)=(\sigma\cap I)\sqcup (\sigma\cap J)$.

The (simplicial) {\it Alexander dual} $K^*_V$ for a simplicial complex $K$ on an index set $V$ such that 
$K\neq\Delta^V$ is defined by
$
K^*_V=\{\sigma\subset V\ | \ \sigma^c\not\in K\}.
$
 If $V=V(K)$ or $V$ is clear from the context, we simply write $K^*$ for $K^*_V$. 
It is easy to see that for a subset $I$ of $V$,
$I$ is a facet of $K$ if and only if $I^c$ is a minimal-non face of $K^*$. The restriction
of $K^*$ and the dual to a link of $K$ are related by the following formula:
$(K^*)_I=(\link_K(I^c))^*.$

Next we review polyhedral products, which are a generalization of moment-angle complexes.

\bigskip

{\bf Definition 2.1.} Let $K$ be a simplicial complex on $[m]$, and $(\X,\A)$ be a sequence of pairs of
based spaces $\{(X_i,A_i)\}_{i\in[m]}$. The polyhedral product $Z_K(\X,\A)$ is defined by
\[
Z_K(\X,\A)=\bigcup_{\sigma\in K}(\X,\A)^\sigma\quad (\subset X_1\times\cdots\times X_m),
\]
where $(\X,\A)^\sigma=Y_1\times\cdots\times Y_m$, with $Y_i=X_i$ for $i\in\sigma$ and $A_i$
for $i\not\in\sigma$. If $(X_i,A_i)=(X,A)$ for all $i\in[m]$, then we write $Z_K(X,A)$ for $Z_K(\X,\A)$.

\bigskip

The polyhedral product $Z_K(D^2,S^1)$ is the moment-angle complex of $K$, and is written simply as $Z_K$.
We refer the reader to \cite{IK14, BP} for further examples of polyhedral products.
In this paper, we are interested in the
homotopy types of polyhedral products $Z_K(C\X,\X)$.

If $K$ has a ghost vertex, then 
\[
Z_K(\X,\A)=Z_{K_{V(K)}}(\X_{V(K)},\A_{V(K)})\times \prod_{i\in V(K)^c}A_i,
\]
where $\X_I$ is the sub-sequence of $\X$ indexed by $I\subset[m]$. In particular, the moment-angle 
complex $Z_K=Z_{K_{V(K)}}\times (S^1)^{|V(K)^c|}$. If $|V(K)^c|>1$ or $Z_{K_{V(K)}}$ is 
not contractible, $K$ is non-Golod.  When we consider Golodness of simplicial complexes, therefore, 
we usually restrict ourselves to simplicial complexes without ghost vertices. 
Nevertheless, in this paper we allow a simplicial complex to have a ghost vertex. The reason for this is 
that we mainly consider dual simplicial complexes as in Theorems 2.3, 2.5 and 2.7. 
These dual simplicial complexes do not have ghost vertices even if the initial simplicial complexes have.  

\bigskip 

{\bf Example 2.2.} In \cite{P65}, Porter proved the following homotopy equivalence:
\begin{eqnarray*}
Z_{\partial\Delta^{[m]}}(C\X,\X)&=&\bigcup_{i=1}^mCX_1\times\cdots\times
 CX_{i-1}\times X_i\times CX_{i+1}\times\cdots\times CX_m\\
&\simeq&X_1*\cdots*X_m\\
&\simeq&\Sigma^{m-1}X_1\wedge\cdots\wedge X_m.
\end{eqnarray*}
He used this homotopy equivalence to define a higher order Whitehead product.

\bigskip

{\bf Theorem 2.3.} {\it Let $K$ and $L$ be simplicial complexes on disjoint index sets $V$ and $W$, respectively.
If $K\neq\Delta^V$ and $L\neq\Delta^W$, then the BBCG decomposition is desuspendable for $(K\ast L)^*$.
}

\bigskip

{\bf Corollary 2.4.} {\it  Let $K$ and $L$ be simplicial complexes on disjoint index sets $V$ and $W$, respectively.
If $K\neq\Delta^V$ and $L\neq\Delta^W$, then $(K\ast L)^*$ is Golod by (1.1).
}

\bigskip

Because the Stanley-Reisner ideal of  $(K*L)^*$ is the product of those of $K^*$ and $L^*$,
this corollary provides a topological proof of a classical result given in \cite{HS}. 
The reader may also refer to \cite{D} and \cite{Ka}.
We also remark that $(\Delta^V*L)^*$ is Golod if and only if
$L^*$ is Golod, because $(\Delta^V*L)^*=\Delta^V*L^*$ (see Lemma 3.1).

\bigskip

{\bf Theorem 2.5.} {\it Let $K$ and $L$ be simplicial complexes of non-negative dimension
on index sets $V$ and $W$, respectively.
Assume that $\alpha=V\cap W$ is a common face of $K$ and $L$.
If $K\neq\Delta^V$ or $L\neq \Delta^W$, and $\alpha$ is neither a facet of $K$ nor $L$, then
the BBCG decomposition is desuspendable for $(K\cup_{\alpha}L)^*$.
}

\bigskip

{\bf Corollary 2.6.} {\it Let $K$ and $L$ be simplicial complexes that satisfy the
same conditions as stated in Theorem 2.5. If $\alpha$ is neither a facet of $K$ nor $L$, 
then $(K\cup_\alpha L)^*$ is Golod by (1.1).
}

\bigskip

It is natural to ask whether $(K\cup_\alpha L)^*$ is still Golod even when $\alpha$ is a facet of $K$ or $L$
in Corollary 2.6. In general, the answer is that this does not hold. 
In fact, we can construct many non-Golod simplicial complexes $(\Delta^V\cup_\alpha L)^*$,
such as in Example 5.3.4. 
A necessary and sufficient condition for $(K\cup_{\alpha}L)^*$ being non-Golod is given by Corollary 5.3.3. 
Incidentally, $(\Delta^V\cup_\alpha\Delta^W)^*$ is a non-Golod simplicial complex if $\alpha\neq V$ and 
$\alpha\neq W$, since $(\Delta^V\cup_\alpha\Delta^W)^*=\partial\Delta^{V-\alpha}*\Delta^\alpha*\partial\Delta^{W-\alpha}$
(see Lemma 4.1).

To construct a simplicial complex satisfying Theorem 1.3, we first need to fix some notation.

Let $K$ be a simplicial complex on an index set $V$, with facets $F_1,\cdots,F_k$. We take new vertices
$v_1,\cdots,v_k$, and define a new simplicial complex $F(K)$ on the index set $V\sqcup\{v_1,\cdots,v_k\}$
with facets $F_1+v_1,\cdots,F_k+v_k$, where $F_i+v_i=F_i\sqcup\{v_i\}$ is the disjoint union of two sets. 
Then, $K$ is a subcomplex of $F(K)$ and $|K|$ is a deformation retract of $|F(K)|$.

For two simplicial complexes $K$ and $L$, we define a ``product'' $K\boxtimes L$ as follows. Define a
liner order $\leq$ on the vertex sets of $K$ and $L$.
The vertex set of $K\boxtimes L$ is $V(K)\times V(L)$. An $n$-simplex is a set
$\{(x_0,y_0),\cdots,(x_n,y_n)\}$ such that $x_0\leq\cdots\leq x_n$, $y_0\leq\cdots\leq y_n$,
$\{x_0,\cdots,x_n\}$ is a simplex of $K$, and $\{y_0,\cdots,y_n\}$ is a simplex of $L$. It is well-known
that $|K\boxtimes L|$ is homeomorphic to $|K|\times|L|$. If $v$ is a vertex of $L$,
then the subcomplex $K\boxtimes\Delta^{\{v\}}$ of $K\boxtimes L$ is abbreviated as $K\boxtimes v$.

By $S^n_k$ we denote a triangulation of an $n$-sphere $S^n$ with $k$-vertices.
It follows from the simplicial approximation theorem that there is a simplicial map $\eta_k:S^3_k\to S^2_4$
for sufficiently large $k$ whose geometrical realization $|\eta_k|:|S^3_k|\to |S^2_4|$ is homotopic 
to the Hopf map $\eta:S^3\to S^2$. In fact, we can choose $k=12$ in this case, by \cite{MS}.

We consider the simplicial set $\Delta^1$ as the full simplex on $[2]$.
By $S^3_k\boxtimes \Delta^1\cup_{\eta_k} S^2_4$
we denote the simplicial complex obtained from the disjoint union of two simplicial complexes
$S^3_k\boxtimes \Delta^1$ and $S^2_4$, given by identifying $(v,2)\in V(S^3_k)\times [2]$ with
$\eta_k(v)\in V(S^2_4)$.
We embed $S^3_k$ into $S^3_k\boxtimes\Delta^1\cup_{\eta_k} S^2_4$ by applying the map $v\mapsto (v,1)$.

We set $V$ to be the vertex set of the union of two simplicial complexes 
$S^3_k\boxtimes \Delta^1\cup_{\eta_k} S^2_4$ and $F(S^3_k)$ 
along $S^3_k$.  
Finally, we take new vertices $v_0,w_1,w_2$ and set
\[
K=\Delta^{V+v_0}\cup ((S^3_k\boxtimes
\Delta^1\cup_{\eta_k} S^2_4)\cup F(S^3_k))*\Delta^{\{w_1\}}
\cup S^3_k*\Delta^{\{w_1,w_2\}},
\]
which is the union of $\Delta^{V+v_0}$ and
$\Delta^V\cup((S^3_k\boxtimes \Delta^1\cup_{\eta_k} S^2_4)\cup F(S^3_k))*\Delta^{\{w_1\}}
\cup S^3_k*\Delta^{\{w_1,w_2\}}$.

\bigskip

{\bf Theorem 2.7.} {\it $K^*$ is a Golod simplicial complex that is not stably
homotopy Golod. Moreover, $Z_{K^*}$ is torsion free if $K$ is constructed from the map $\eta_{12}:S^3_{12}
\to S^2_4$ defined in \cite{MS}.}

\section{Proof of Theorem 2.3.}

In this section, we prove Theorem 2.3. We begin by stating some elementary lemmas, for which the proofs are omitted.

\bigskip

{\bf Lemma 3.1.} {\it Let $K$ and $L$ be simplicial complexes with disjoint index sets $V$ and $W$, respectively.
Then, $(K\ast L)^*_{V\sqcup W}=K^*_V\ast\Delta^W\cup\Delta^V\ast L^*_W.$
}

\bigskip

{\bf Lemma 3.2.} {\it Let $K$ and $L$ be simplicial complexes with disjoint index sets, 
and let $K_i$ and $L_i$ for $i=1,2$ be subcomplexes of $K$ and $L$, respectively. 
Then, $$(K_1\ast L_1)\cap(K_2\ast L_2)=(K_1\cap K_2)\ast(L_1\cap L_2).$$
}

\bigskip

{\bf Lemma 3.3.} {\it Let $K$ be a simplicial complex with two subcomplexes $K_1$ and $K_2$.
If \\ $K=K_1\cup K_2$, then 
\[Z_K(C\X,\X)=Z_{K_1}(C\X,\X)\cup Z_{K_2}(C\X,\X)
\]
and
\[
Z_{K_1}(C\X,\X)\cap Z_{K_2}(C\X,\X)=Z_{K_1\cap K_2}(C\X,\X).
\]
}

\bigskip

{\bf Lemma 3.4.} {\it Let $K$ and $L$ be simplicial complexes with disjoint index sets $V$ and $W$, respectively.  Then,
\[
Z_{K*L}(C\X,\X)=Z_K(C\X_V,\X_V)\times  Z_L(C\X_W,\X_W).
\]
}

\bigskip

{\bf Proposition 3.5.} {\it Let $K$ be a simplicial complex on $[m]$ and $\X$ be a sequence of based
CW-complexes. If $Z_K(C\X,\X)$ is a simply connected co-H-space,
then 
\[
Z_K(C\X,\X)\simeq \bigvee_{I\subset[m]}\Sigma|K_I|\wedge\widehat{X}^I.
\]
}

{\bf Proof.} First, we show that $\bigvee_{I\subset[m]}\Sigma|K_I|\wedge\widehat{X}^I$ is also simply connected.
By Theorem 1.2 we have that $H_1( \bigvee_{I\subset[m]}\Sigma|K_I|\wedge\widehat{X}^I)\cong H_1(Z_K(C\X,\X))=0$.
Because $\bigvee_{I\subset[m]}\Sigma|K_I|\wedge\widehat{X}^I$ is a suspension, its fundamental group
is a free group. Thus, $\pi_1(\bigvee_{I\subset[m]}\Sigma|K_I|\wedge\widehat{X}^I)=0$.

For a subset $I\subset[m]$, the canonical projection $p_I:\prod_{i\in[m]}CX_i\to \prod_{i\in I}CX_i$
induces a map $Z_K(C\X,\X)\to Z_{K_I}(C\X_I,\X_I)$, which is also denoted by $p_I$. Let $I=\{i_1,\cdots,i_k\}\subset[m]$. We define a subset
$Z_{K_I}(C\X_I,\X_I)^\prime$ of $Z_{K_I}(C\X_I,\X_I)$ by the equation
\[
Z_{K_I}(C\X_I,\X_I)^\prime=\{(x_{i_1},\cdots,x_{i_k})\in Z_{K_I}(C\X_I,\X_I) \ |\ \text{$x_{i_j}=*$ for some $j$}\}. 
\]
In \cite{IK14}, it is shown that $Z_{K_I}(C\X_I,\X_I)/Z_{K_I}(C\X_I,\X_I)^\prime\simeq
\Sigma|K_I|\wedge\widehat{X}^I$. Now, we consider the composite of maps
\begin{multline*}
f:Z_K(C\X,\X)\to \bigvee^{2^m}Z_K(C\X,\X)\xrightarrow{\vee_{I\subset[m]} p_I}\bigvee_{I\subset[m]}
Z_{K_I}(C\X_I,\X_I)\\
\to \bigvee_{I\subset[m]}Z_{K_I}(C\X_I,\X_I)/Z_{K_I}(C\X_I,\X_I)^\prime
\simeq \bigvee_{I\subset[m]}\Sigma|K_I|\wedge\widehat{X}^I,
\end{multline*}
where the first map is the iterated co-multiplication of $Z_K(C\X,\X)$. Because $Z_K(C\X,\X)$ and
$ \bigvee_{I\subset[m]}\Sigma|K_I|\wedge\widehat{X}^I$ are simply connected CW-complexes, to prove that $f$ is homotopy equivalent 
it suffices to show that $f$ induces a homology isomorphism.
In \cite{IK15}, it is shown that $\Sigma f$
is a homotopy equivalence. In particular, $f$ induces a homology isomorphism, and thus we complete the
proof. \hfill $\Box$

\bigskip

{\bf Proof of Theorem 2.3.} By Lemmas 3.1 and 3.2, we have that
$(K\ast L)^*=K^*\ast\Delta^W\cup\Delta^V\ast L^*$ and $K^*\ast\Delta^{V(L)}\cap\Delta^{V(K)}\ast L^*=K^*\ast L^*$.
Therefore, from Lemma 3.3 we obtain the following push-out diagram of spaces:
\[
\begin{CD}
Z_{K^*\ast L^*}(C\X,\X)		@>>>	Z_{K^*\ast\Delta^W}(C\X,\X)\\
@VVV						@VVV\\
Z_{{\Delta^V}\ast L^*}(C\X,\X)	@>>>	Z_{(K\ast L)^*}(C\X,\X).
\end{CD}
\]
Here, we remark that $K^*\ast L^*$ and $K^*\ast\Delta^W$ are non-void simplicial complexes,
because we assume that $K\neq\Delta^V$ and $L\neq\Delta^W$.
By Lemma 3.4, the above push-out diagram is equivalent to the following push-out diagram:
\[
\begin{CD}
Z_{K^*}(C\X_V,\X_V)\times Z_{L^*}(C\X_W,\X_W)	@>>>	
Z_{K^*}(C\X_V,\X_V)\times\prod_{w\in W}CX_w\\
@VVV				@VVV\\
\prod_{v\in V}CX_v\times Z_{L^*}(C\X_W,\X_W)	@>>>
	Z_{(K\ast L)^*}(C\underline{X},\underline{X}).
\end{CD}.
\]
Because $\prod_{v\in V}CX_v$ and $\prod_{w\in W}CX_w$ are contractible, the above diagram yields the following
homotopy equivalences:
\begin{eqnarray*}
Z_{(K\ast L)^*}(C\X,\X)&\simeq& Z_{K^*}(C\X_V,\X_V)\ast Z_{L^*}(C\X_W,\X_W)\\
&\simeq&\Sigma Z_{K^*}(C\X_V,\X_V)\wedge Z_{L^*}(C\X_W,\X_W).\\
\end{eqnarray*}
By Theorem 1.2, $\Sigma Z_{K^*}(C\X_V,\X_V)$ is a double suspension, which implies that
$Z_{(K\ast L)^*}(C\X,\X)$ is also a double suspension. By invoking Proposition 3.5,
we complete the proof. \hfill $\Box$

\section{Proof of Theorem 2.5.}

In this section, we prove Theorem 2.5. Again, we begin by stating some elementary lemmas.

\bigskip

{\bf Lemma 4.1.} {\it Let $K$ and $L$ be simplicial complexes without ghost vertices 
on index sets $V$ and $W$, respectively.
Let $\alpha=V\cap W$ be a common face of $K$ and $L$.
Then, 
\[
(K\cup_\alpha L)^*
=(\partial\Delta^{V-\alpha}\ast\Delta^{\alpha}\ast\partial\Delta^{W-\alpha})
\cup(K^*\ast\Delta^{W-\alpha})\cup(\Delta^{V-\alpha}\ast L^*).
\]
}

\bigskip

{\bf Proof.} For $u\in V-\alpha$ and $v\in W-\alpha$, $\{u,v\}$ is a minimal non-face of $K\cup_\alpha L$, 
and a minimal non-face of $K$ or $L$ is also a minimal non-face of $K\cup_\alpha L$.  
For any $u\in V-\alpha$ and $v\in W-\alpha$,
$V\cup W-\{u,v\}$ is a facet of $(K\cup_\alpha L)^*$. Furthermore, for any minimal non-face
$\sigma$ of $K$ or $L$, $V\cup W- \sigma$ is a facet of $(K\cup_\alpha L)^*$.
This implies the desired equality of the two complexes. \hfill $\Box$

\bigskip

{\bf Lemma 4.2.} {\it If $K$ is a simplicial complex with a ghost vertex, then the BBCG decomposition
for $K^*$ is desuspendable. In particular, $K^*$ is Golod.}

\bigskip

{\bf Proof.} Let $K$ be a simplicial complex on $[m]$ and $v$ be a ghost vertex of $K$.
Then, $[m]-v$ is a facet of $K^*$, and thus $\dim K^*\geq m-2$.  Then, it follows from Theorem 1.2 and
Proposition 3.5 of \cite{IK14} that the BBCG decomposition for $K^*$ is desuspendable. \hfill $\Box$

\bigskip

{\bf Lemma 4.3.} {\it Let $\alpha$ be a face of a simplicial complex $K$ on an index set $V$.
If $K\neq\Delta^V$ and $\alpha$ is not a facet of $K$, then the inclusion map
\[
Z_{K^*}(C\X,\X)\to Z_{(\Delta^\alpha)^*}(C\X,\X)=
Z_{\partial\Delta^{V-\alpha}\ast\Delta^\alpha}(C\X,\X)
\]
is null homotopic.
}

\bigskip

{\bf Proof.} Because $\alpha$ is not a facet of $K$, there is a face $\beta$ of $K$
such that $\alpha\subsetneq\beta$. Then, \\$\Delta^\alpha\subsetneq\Delta^\beta\subset K$, 
which implies that $K^*\subset(\Delta^\beta)^*\subsetneq(\Delta^\alpha)^*$. That is, \\
$
K^*\subset \partial\Delta^{V-\beta}\ast \Delta^\beta\subsetneq\partial\Delta^{V-\alpha}\ast\Delta^\alpha.
$
Therefore, the inclusion $Z_{K^*}(C\X,\X)\hookrightarrow
Z_{\partial\Delta^{V-\alpha}\ast\Delta^\alpha}(C\X,\X)$ factors
as $Z_{K^*}(C\X,\X)\rightarrow Z_{\partial\Delta^{V-\beta}\ast\Delta^\beta}(C\X,\X)
\rightarrow Z_{\partial\Delta^{V-\alpha}\ast\Delta^\alpha}(C\X,\X)$. To show that the inclusion
$Z_{K^*}(C\X,\X)\rightarrow Z_{\partial\Delta^{V-\alpha}\ast\Delta^\alpha}(C\X,\X)$ is
null homotopic, it is sufficient to show that \\
$Z_{\partial\Delta^{V-\beta}\ast\Delta^\beta}(C\X,\X)\rightarrow Z_{\partial\Delta^{V-\alpha}\ast\Delta^\alpha}(C\X,\X)$ is null homotopic.
We have the following homotopy commutative diagram  
\[
\begin{CD}
 Z_{\partial\Delta^{V-\beta}\ast\Delta^\beta}(C\X,\X)	=
Z_{\partial\Delta^{V-\beta}}(C\X_{V-\beta},\X_{V-\beta})\times \prod_{j\in\beta}CX_j
@>\simeq>> *_{i\in V-\beta}X_i,\\
@VVV                       @VVV   \\
 Z_{\partial\Delta^{V-\alpha}\ast\Delta^\alpha}(C\X,\X)	=
Z_{\partial\Delta^{V-\alpha}}(C\X_{V-\alpha},\X_{V-\alpha})\times \prod_{j\in\alpha}CX_j
@>\simeq>> (*_{i\in V-\beta}X_i)*(*_{i\in\beta-\alpha}X_i),\\
\end{CD}
\]
where the horizontal maps are homotopy equivalences by Example 2.2 and the right vertical map 
is the canonical inclusion $X\to X*Y$.  Since this map is null homotopic, \\
$Z_{\partial\Delta^{V-\beta}\ast\Delta^\beta}(C\X,\X)\rightarrow Z_{\partial\Delta^{V-\alpha}\ast\Delta^\alpha}(C\X,\X)$ is also null homotopic, 
and thus we complete the proof. \hfill $\Box$

\bigskip

In addition to the above, we require the following lemma to prove Theorem 2.5.

\bigskip

{\bf Lemma 4.4 (Lemma 3.2 of \cite{IK13} with $B=*$).} {\it Define $Q$ as the push-out
\[
\begin{CD}
A\times C @>\iota\times1>>CA\times C\\
@VV1\times*V @VVV\\
A\times D@>j>>      Q,\\
\end{CD}
\]
where $\iota:A\to CA$ is the inclusion. Then, the homotopy equivalence
\[
Q\xrightarrow{\simeq} (A\ltimes D)\vee\Sigma(A\wedge C)
\]
holds, which is natural with respect to $A,\ C,$ and $D$, where $X\ltimes Y=(X\times Y)/(X\times *)$.
Moreover, the inclusion map $j: A\times D\to Q$ is homotopic to the following composite of maps 
\[
 A\times D\to  A\ltimes D\to  (A\ltimes D)\vee\Sigma(A\wedge C)\simeq Q,
\]
where the first map is the collapsing map and the second map is the inclusion. 
}

\bigskip

{\bf Proof of Theorem 2.5.} If $K$ or $L$ has a ghost vertex, then $K\cup_\alpha L$ also has a ghost vertex.
In this case, the BBCG decomposition is desuspendable, by Lemma 4.2.
Therefore, we assume that $K$ and $L$ do not have any ghost vertices. In the following proof, $Z_K(C\X,\X)$ is
abbreviated as $Z_K$.

First, we will show that if $K=\Delta^V$ and $L\neq \Delta^W$, then
\[
Z_{(\Delta^V\cup_\alpha L)^*}	\simeq	(Z_{\partial\Delta^{V-\alpha}}\ltimes Z_{\partial\Delta^{W-\alpha}})\vee\Sigma(Z_{\partial\Delta^{V-\alpha}}\wedge Z_{L^*}).
\]
Here we remark that we need the assumption that $L\neq \Delta^W$ and that $\alpha$ is not a facet of $L$ 
to apply Lemma 4.3. 

It follows from Lemma 4.1 that
$
(\Delta^V\cup_\alpha L)^*=
(\partial\Delta^{V-\alpha}\ast\Delta^\alpha\ast\partial\Delta^{W-\alpha})\cup(\Delta^{V-\alpha}\ast L^*).
$
Then, by Lemma 3.3 we have the push-out diagram of spaces
\[
\xymatrix{
Z_{\partial\Delta^{V-\alpha}\ast L^*}\ar[r]\ar[d]		
&	Z_{\Delta^{V-\alpha}\ast L^*}\ar[d]\\
Z_{\partial\Delta^{V-\alpha}\ast\Delta^\alpha\ast\partial\Delta^{W-\alpha}}\ar[r]	&	Z_{(\Delta^V\cup_\alpha L)^*},
}
\]
which by Lemma 3.4 is equivalent to the following push-out diagram:
\[
\begin{CD}
Z_{\partial\Delta^{V-\alpha}}\times Z_{L^*}				@>>>	Z_{\Delta^{V-\alpha}}\times Z_{L^*}\\
@V{\mathrm id}\times\mathrm{incl}VV					@VVV\\
Z_{\partial\Delta^{V-\alpha}}\times Z_{\Delta^\alpha}\times Z_{\partial\Delta^{W-\alpha}}								@>>>Z_{(\Delta^V\cup_\alpha L)^*}.
\end{CD}
\]
Since we assumed that $L\neq \Delta^W$ and that $\alpha$ is not a facet of $L$, 
by Lemma 4.3 the inclusion map 
$Z_{L^*}\hookrightarrow Z_{\Delta^\alpha}\times Z_{\partial\Delta^{W-\alpha}}$
is null-homotopic. Therefore, $Z_{(\Delta^V\cup_\alpha L)^*}$ is homotopy equivalent to
the push-out $P$ of the following diagram:
\begin{equation}
\begin{CD}
Z_{\partial\Delta^{V-\alpha}}\times Z_{L^*}				@>>>	Z_{\Delta^{V-\alpha}}\times Z_{L^*}\\
@V{\mathrm id}\times\ast VV															@VVV\\
Z_{\partial\Delta^{V-\alpha}}\times Z_{\Delta^\alpha}\times Z_{\partial\Delta^{W-\alpha}}								@>j>>	P. 
\end{CD}
\end{equation}
Since $Z_{\Delta^{V-\alpha}}=CZ_{\partial\Delta^{V-\alpha}}$, by Lemma 4.4 
$P$ is homotopy equivalent to
\[
(Z_{\partial\Delta^{V-\alpha}}\ltimes(Z_{\Delta^\alpha}\times  Z_{\partial\Delta^{W-\alpha}}))\vee
\Sigma(Z_{\partial\Delta^{V-\alpha}}\wedge Z_{L^*})\simeq
(Z_{\partial\Delta^{V-\alpha}}\ltimes
Z_{\partial\Delta^{W-\alpha}})\vee\Sigma(Z_{\partial\Delta^{V-\alpha}}\wedge Z_{L^*}),
\]
and $j$ in the diagram (4.1) can be identified with the following composite of canonical maps:
\begin{multline}
Z_{\partial\Delta^{V-\alpha}}\times Z_{\Delta^\alpha}\times Z_{\partial\Delta^{W-\alpha}}	\to
Z_{\partial\Delta^{V-\alpha}}\times Z_{\partial\Delta^{W-\alpha}}	\to\\ 
Z_{\partial\Delta^{V-\alpha}}\ltimes Z_{\partial\Delta^{W-\alpha}}
\hookrightarrow	
(Z_{\partial\Delta^{V-\alpha}}\ltimes Z_{\partial\Delta^{W-\alpha}})\vee\Sigma(Z_{\partial\Delta^{V-\alpha}}\wedge Z_{L^*}).
\end{multline}
Thus, the following homotopy equivalence holds:
\[
Z_{(\Delta^V\cup_\alpha L)^*}\simeq(Z_{\partial\Delta^{V-\alpha}}\ltimes Z_{\partial\Delta^{W-\alpha}})\vee\Sigma(Z_{\partial\Delta^{V-\alpha}}\wedge Z_{L^*}).
\]
This homotopy equivalence induces the following homotopy equivalences:
\begin{eqnarray*}
Z_{(\Delta^V\cup_\alpha L)^*}&\simeq&
(Z_{\partial\Delta^{V-\alpha}}\ltimes Z_{\partial\Delta^{W-\alpha}})\vee\Sigma(Z_{\partial\Delta^{V-\alpha}}\wedge Z_{L^*})\\
&\simeq&Z_{\partial\Delta^{W-\alpha}}\vee (Z_{\partial\Delta^{V-\alpha}}\wedge
Z_{\partial\Delta^{W-\alpha}})\vee \Sigma(Z_{\partial\Delta^{V-\alpha}}\wedge Z_{L^*})\\
&\simeq& Z_{\partial\Delta^{W-\alpha}}\vee (Z_{\partial\Delta^{V-\alpha}}\wedge
Z_{\partial\Delta^{W-\alpha}})\vee \bigvee_{J\subset W}\Sigma Z_{\partial\Delta^{V-\alpha}}\wedge
\Sigma |(L^*)_J|\wedge\widehat{X}^J\\
&\simeq& Z_{\partial\Delta^{W-\alpha}}\vee (Z_{\partial\Delta^{V-\alpha}}\wedge
Z_{\partial\Delta^{W-\alpha}})\vee \bigvee_{J\subset W}\Sigma^3|\partial\Delta^{V-\alpha}|
\wedge|(L^*)_J|\wedge\widehat{X}^{V-\alpha+J}\\
&\simeq& \bigvee_{I\subset[m]}\Sigma|((\Delta^V\cup_\alpha L)^*)_I|
\wedge \widehat{X}^I,\\
\end{eqnarray*}
where in the third homotopy equivalence we used the BBCG decomposition for $\Sigma Z_{L^*}$. 
To show that the last homotopy equivalence holds we need the homotopy type of 
$|((\Delta^V\cup_\alpha L)^*)_I|$ for 
$I\subset[m]$. Since $(\Delta^V\cup_\alpha L)^*=
(\partial\Delta^{V-\alpha}\ast\Delta^\alpha\ast\partial\Delta^{W-\alpha})\cup(\Delta^{V-\alpha}\ast L^*),$
we have the following push-out diagram for $I=J^\prime\sqcup J$ with $J^\prime\subset V-\alpha$ and 
$J\subset W$: 
\[
\begin{CD}
|(\partial\Delta^{V-\alpha})_{J^\prime}*(L^*)_J| @>>> |(\Delta^{V-\alpha})_{J^\prime}*(L^*)_J|\\
@VVV                   @VVV             \\
|(\partial\Delta^{V-\alpha})_{J^\prime}*(\Delta^\alpha*\partial\Delta^{W-\alpha})_J|@>>> 
|((\Delta^V\cup_\alpha L)^*)_I|.
\end{CD}
\]
Since $\alpha$ is not a facet of $L$, there is a face $\beta$ of $L$ and we have inclusions $\Delta^\alpha\subsetneq\Delta^\beta\subset L$. This implies that $|(L^*)_J|\subset 
|(\Delta^\beta*\partial\Delta^{W-\beta})_J|\subset |(\Delta^\alpha*\partial\Delta^{W-\alpha})_J|$. Since 
the inclusion map $|(\Delta^\beta*\partial\Delta^{W-\beta})_J|\hookrightarrow
|(\Delta^\alpha*\partial\Delta^{W-\alpha})_J|$ is null homotopic, the left vertical map of the diagram above 
is null homotopic. Thus we have 
\[
 |((\Delta^V\cup_\alpha L)^*)_I|\simeq 
|(\partial\Delta^{V-\alpha})_{J^\prime}*(\Delta^\alpha*\partial\Delta^{W-\alpha})_J|
\vee  |(\Delta^{V-\alpha})_{J^\prime}*(L^*)_J|/|(\partial\Delta^{V-\alpha})_{J^\prime}*(L^*)_J|. 
\]
If $J^\prime\neq\emptyset,V-\alpha$, then $|(\partial\Delta^{V-\alpha})_{J^\prime}|$ is contractible, and  
$|((\Delta^V\cup_\alpha L)^*)_I|$ is also contractible. If $J^\prime=\emptyset$, then 
$|((\Delta^V\cup_\alpha L)^*)_I|=|(\Delta^\alpha*\partial\Delta^{W-\alpha})_J|$ which is 
$|\partial\Delta^{W-\alpha}|$ if $J=W-\alpha$, and is contractible otherwise. Finally we obtain the 
following homotopy equivalence, where $|\partial\Delta^{V-\alpha}*(\partial\Delta^{W-\alpha})_J|$ 
is contractible unless $J=W-\alpha$. 
\begin{multline*}
|((\Delta^V\cup_\alpha L)^*)_I|\simeq\\
\begin{cases}
|\partial\Delta^{W-\alpha}| &\text{for $I=W-\alpha$},\\
|\partial\Delta^{V-\alpha}*(\partial\Delta^{W-\alpha})_J|\vee\Sigma|\partial\Delta^{V-\alpha}*(L^*)_J|
&\text{for $I=(V-\alpha)\sqcup J$ for some $J\subset W$},\\
* &\text{otherwise}.
\end{cases}
\end{multline*}
Here, we remark that we can apply Proposition 3.5 if $|W-\alpha|\geq3$.

Similarly, if $K\neq\Delta^{V}\text{ and }L= \Delta^{W}$ then we have
\[
Z_{(K\cup_\alpha \Delta^W)^*}\simeq(Z_{\partial\Delta^{V-\alpha}}\rtimes Z_{\partial\Delta^{W-\alpha}})\vee\Sigma(Z_{K^*}\wedge Z_{\partial\Delta^{W-\alpha}})
\simeq \bigvee_{I\subset[m]}\Sigma|((K\cup_\alpha \Delta^W)^*)_I|\wedge 
\widehat{X}^I.
\]

Next, we consider the case that $K\neq\Delta^{V}\text{ and }L\neq\Delta^{W}$.
Then, we have the push-out diagram
\[
\begin{CD}
(\partial\Delta^{V-\alpha}\ast\Delta^\alpha\ast\partial\Delta^{W-\alpha})		@>>>	 (\partial\Delta^{V-\alpha}\ast\Delta^\alpha\ast\partial\Delta^{W-\alpha})\cup (\Delta^{V-\alpha}\ast L^*)\\
@VVV						@VVV\\
(\partial\Delta^{V-\alpha}\ast\Delta^\alpha\ast\partial\Delta^{W-\alpha})\cup (K^*\ast\Delta^{W-\alpha})	@>>>	(K\cup_\alpha L)^*,
\end{CD}
\]
which induces the following push-out diagram of spaces by Lemmas 3.3 and 3.4:
\[
\begin{CD}
Z_{\partial\Delta^{V-\alpha}}\times Z_{\Delta^\alpha}\times Z_{\partial\Delta^{W-\alpha}}							@>>>	 Z_{(\partial\Delta^{V-\alpha}\ast\Delta^\alpha\ast\partial\Delta^{W-\alpha})\cup (\Delta^{V-\alpha}\ast L^*)}\\
@VVV							@VVV\\
Z_{(\partial\Delta^{V-\alpha}\ast\Delta^\alpha\ast\partial\Delta^{W-\alpha})\cup (K^*\ast\Delta^{W-\alpha})}			@>>>	Z_{(K\cup_\alpha L)^*}.
\end{CD}.
\]
Because the top arrow of the diagram above is a closed cofibration, 
$Z_{(K\cup_\alpha L)^*}$ is homotopy equivalent to the homotopy push-out of the following diagram:
\[
\begin{CD}
Z_{\partial\Delta^{V-\alpha}}\times Z_{\Delta^\alpha}\times Z_{\partial\Delta^{W-\alpha}}							@>>>	 Z_{(\partial\Delta^{V-\alpha}\ast\Delta^\alpha\ast\partial\Delta^{W-\alpha})\cup (\Delta^{V-\alpha}\ast L^*)}\\
@VVV							@.\\
Z_{(\partial\Delta^{V-\alpha}\ast\Delta^\alpha\ast\partial\Delta^{W-\alpha})\cup (K^*\ast\Delta^{W-\alpha})}			@.
\end{CD}.
\]
Since we have the following homotopy equivalences
\begin{eqnarray*}
Z_{(\partial\Delta^{V-\alpha}\ast\Delta^\alpha\ast\partial\Delta^{W-\alpha})\cup (\Delta^{V-\alpha}\ast L^*)}&\simeq&	
(Z_{\partial\Delta^{V-\alpha}}\ltimes (Z_{\Delta^\alpha}\times Z_{\partial\Delta^{W-\alpha}}))\vee\Sigma (Z_{\partial\Delta^{V-\alpha}}\wedge Z_{L^*})\\
&\simeq&(Z_{\partial\Delta^{V-\alpha}}\ltimes Z_{\partial\Delta^{W-\alpha}})\vee\Sigma (Z_{\partial\Delta^{V-\alpha}}\wedge Z_{L^*}),\\
Z_{(\partial\Delta^{V-\alpha}\ast\Delta^\alpha\ast\partial\Delta^{W-\alpha})\cup (K^*\ast\Delta^{W-\alpha})}	&\simeq&	
((Z_{\partial\Delta^{V-\alpha}}\times Z_{\Delta^\alpha})\rtimes Z_{\partial\Delta^{W-\alpha}})\vee \Sigma (Z_{K^*}\wedge
Z_{\partial\Delta^{W-\alpha}})\\
&\simeq&(Z_{\partial\Delta^{V-\alpha}}\rtimes Z_{\partial\Delta^{W-\alpha}})\vee \Sigma (Z_{K^*}\wedge Z_{\partial\Delta^{W-\alpha}}),
\end{eqnarray*}
and homotopy push-out preserves homotopy equivalences,  
$Z_{(K\cup_\alpha L)^*}$ is homotopy equivalent to the homotopy push-out of the following diagram:
\[
\begin{CD}
Z_{\partial\Delta^{V-\alpha}}\times Z_{\Delta^\alpha}\times Z_{\partial\Delta^{W-\alpha}}							@>>j_1>	 (Z_{\partial\Delta^{V-\alpha}}\ltimes Z_{\partial\Delta^{W-\alpha}})\vee\Sigma (Z_{\partial\Delta^{V-\alpha}}\wedge Z_{L^*})\\
@VVj_2V\\
(Z_{\partial\Delta^{V-\alpha}}\rtimes Z_{\partial\Delta^{W-\alpha}})\vee \Sigma (Z_{K^*}\wedge Z_{\partial\Delta^{W-\alpha}})
\end{CD}
\]
where $j_1$ and $j_2$ are as in (4.2). It follows that 
\[
Z_{(K\cup_\alpha L)^*}\simeq Q\vee \Sigma (Z_{\partial\Delta^{V-\alpha}}\wedge Z_{L^*})\vee
\Sigma (Z_{K^*}\wedge Z_{\partial\Delta^{W-\alpha}}),
\]
where $Q$ is the homotopy push-out of the following diagram:
\[
\begin{CD}
Z_{\partial\Delta^{V-\alpha}}\times Z_{\partial\Delta^{W-\alpha}}
@>>>Z_{\partial\Delta^{V-\alpha}}\ltimes  Z_{\partial\Delta^{W-\alpha}}\\
@VVV\\
Z_{\partial\Delta^{V-\alpha}}\rtimes Z_{\partial\Delta^{W-\alpha}}
\end{CD}.
\]
where the vertical and horizontal maps are the canonical collapsing maps. Then, it is easy to see that
$Q\simeq Z_{\partial\Delta^{V-\alpha}}\wedge Z_{\partial\Delta^{W-\alpha}}$, and we
have obtained a homotopy equivalence
\begin{equation}
Z_{(K\cup_\alpha L)^*}\simeq
(Z_{\partial\Delta^{V-\alpha}}\wedge Z_{\partial\Delta^{W-\alpha}})\vee
\Sigma (Z_{\partial\Delta^{V-\alpha}}\wedge Z_{L^*})\vee
\Sigma (Z_{K^*}\wedge Z_{\partial\Delta^{W-\alpha}}).
\end{equation}
To complete the proof, we now apply Proposition 3.5. Clearly, (4.3) implies that $Z_{(K\cup_\alpha L)^*}$
is a suspension space. Because $\alpha$ is neither a facet of $K$ nor $L$, we have that 
$|V-\alpha|\geq2$ and $|W-\alpha|\geq2$.
Then, by Example 2.2 we have that $Z_{\partial\Delta^{V-\alpha}}=Z_{\partial\Delta^{V-\alpha}}(C\X,\X)$ 
and \\
   $Z_{\partial\Delta^{W-\alpha}}=Z_{\partial\Delta^{W-\alpha}}(C\X,\X)$ are suspensions.
Thus, (4.3) implies that $Z_{(K\cup_\alpha L)^*}$ is a double suspension, and by Proposition 3.5 
we have the desired homotopy equivalence 
\[
Z_{(K\cup_\alpha L)^*}(C\X,\X)\simeq \bigvee_{I\subset[m]}\Sigma|((K\cup_\alpha L)^*)_I|\wedge 
\widehat{X}^I
\]
and we complete the proof.  \hfill $\Box$

\section{Preliminaries for the proof of Theorem 2.7}

In this section, we review the Alexander duality and elementary collapses of simplicial complexes, and
the Massey products of the Koszul homology.

\subsection{Alexander duality}
First, we review the Alexander duality of simplicial complexes. Although an elementary proof has been
derived \cite{BT}, we follow the classical argument presented in chapter six of Spanier's book \cite{S}.

Let $M$ be a simplicial complex that is a subcomplex of $\partial\Delta^{n+1}$, and
let $L$ be a subcomplex of $M$. We denote the Alexander dual of $L$ and $M$ in
$V(\partial\Delta^{n+1})$ by $L^*,\ M^*$. Then, $M^*$ is a subcomplex of $L^*$, and
the following duality is well-known, where we write $S^n=|\partial\Delta^{n+1}|$:
\begin{multline*}
\gamma:H_q(|L^*|,|M^*|)\cong H_q(|\Sd L^*|,|\Sd M^*|)\cong H_q(|\bar{L}|,|\bar{M}|)\\
\cong H_q(S^n-|L|,S^n-|M|)\cong H^{n-q}(|M|,|L|)
\end{multline*}
where $\Sd K$ denotes the barycentric subdivision of $K$, $\bar{L}$ is the supplement of $L$ in $\partial\Delta^{n+1}$ defined in Definition 2.5.18 of
\cite{M}, and the last isomorphism is the (topological) duality
\[
\gamma_U:H_q(S^n-|L|,S^n-|M|)\cong H^{n-q}(|M|,|L|)
\]
induced by an orientation class $U\in H^n(S^n\times S^n,S^n\times S^n-\Delta(S^n))$ of $S^n$,
where $\Delta(S^n)$ is the diagonal set of $S^n\times S^n$.

Now, we consider the cohomology theory with coefficient $\Z_2=\Z/2\Z$, and the relation between
the duality map $\gamma:H_q(|L^*|,|M^*|;\Z_2)\to H^{n-q}(|M|,|L|;\Z_2)$ and the Steenrod squaring
operations $Sq^i$. For a pair of finite complexes $(X,A)$, the Steenrod squaring operations $Sq^i_*$ on
the homology theory are defined as follows:
\[
<x,Sq^i_*a>=<Sq^ix,a>\quad\text{for $x\in H^q(X,A;\Z_2)$ and $a\in H_{q+i}(X,A;\Z_2)$},
\]
where $<-,->$ denotes the Kronecker product.

\bigskip

{\bf Lemma 5.1.1.} For $a\in H_q(|L^*|,|M^*|)$ and $i>0$, we have that
\[
\gamma(Sq^i_*a)=\sum_{k=0}^{i-1}Sq^{i-k}(\gamma(Sq^k_*a)).
\]

{\bf Proof.} In this proof, we omit the coefficient ring $\Z_2$ in the (co)homology theory. First, we remark that
$Sq^i(U)=0$ for $i>0$. This follows from the fact that the natural restriction map
$H^n(S^n\times S^n,S^n\times S^n-\Delta(S^n))\to \tilde{H}^n(S^n\times S^n)$ is monomorphic.
In the definition of the duality map $\gamma$, all maps except for $\gamma_U$ commute with the squaring
operations. Therefore, it suffices to prove the corresponding formula for $\gamma_U$.

Because $\gamma_U(Sq^i_*a)\in  H^{n-q+i}(|M|,|L|)$, we take an element
$b\in H_{n-q+i}(|M|,|L|)$ and compute the Kronecker product $<\gamma_U(Sq^i_*a),b>$.
Let
\[j:(|M|,|L|)\times(S^n-|L|,S^n-|M|)\to (S^n\times S^n,S^n\times S^n-\Delta(S^n))
\]
 be the inclusion map. Then,
\begin{eqnarray*}
<\gamma_U(Sq^i_*a),b>&=&<j^*(U)/Sq^i_*a,b>\\
&=&<j^*(U),b\times Sq^i_*a>\\
&=&<j^*(U),Sq^i_*(b\times a)-\sum_{k=0}^{i-1}Sq^{i-k}_*b\times Sq^k_*a>\\
&=&<Sq^ij^*(U),b\times a>+\sum_{k=0}^{i-1}<j^*(U),Sq^{i-k}_*b\times Sq^k_*a>\\
&=&\sum_{k=0}^{i-1}<j^*(U)/Sq^k_*a,Sq^{i-k}_*b>\\
&=&<\sum_{k=0}^{i-1}Sq^{i-k}(\gamma_U(Sq^k_*a)),b>,
 \end{eqnarray*}
where in the fifth equation we use the fact that $Sq^ij^*(U)=j^*Sq^i(U)=0$.
Thus, the proof is completed. \hfill $\Box$

\subsection{Elementary collapse}

In this subsection, we briefly review the notion of elementary collapse, which is necessary to prove Theorem 2.7.

\bigskip

{\bf Definition 5.2.1.} A non-empty face $\sigma$ in a simplicial complex $K$ is a {\it free face} if
it is not a facet of $K$ and is contained in exactly one facet of $K$.

An {\it elementary collapse} of $K$ is a simplicial complex $K^\prime$ obtained from $K$ by the removal
of a free face $\sigma$, along with all faces that contain $\sigma$. If there is a sequence of elementary
collapses leading from $K$ to $K^\prime$, we say that $K$ is collapsible onto $K^\prime$, and we use the
notation $K\searrow K^\prime$. Then $|K^\prime|$ is a deformation retract of $|K|$.

\bigskip

{\bf Example 5.2.2.} Let $\Delta^m$ denote the full simplex on the vertex set $[m+1]$. There
exists a sequence of elementary collapses of $\Delta^m\boxtimes\Delta^1$ given by removing
pairs of faces from the top to bottom in the following list:
\begin{eqnarray*}
&&\{\underline{1}\cdots\underline{m+1},\ \underline{1}\cdots\underline{m+1}\overline{m+1}\},\\
&&\cdots\\
&&\{\underline{1}\cdots\underline{i}\overline{i+1}\cdots\overline{m+1},\ \underline{1}\cdots
\underline{i}\overline{i}\ \overline{i+1}\cdots\overline{m+1}\},\\
&&\cdots\\
&&\{\underline{1}\overline{2}\cdots\overline{m+1},\ \underline{1}\overline{1}\ \overline{2}\cdots\overline{m+1}\},\\
\end{eqnarray*}
where $\underline{i}=(i,1)$, $\overline{i}=(i,2)$, and
$\underline{1}\cdots\underline{i}\overline{i+1}\cdots\overline{m+1}$ denotes, for example, the face \\
$\{\underline{1},\cdots,\underline{i},\overline{i+1},\cdots,\overline{m+1}\}$ of $\Delta^m\boxtimes\Delta^1$.
Thus, we see that
$\Delta^m\boxtimes\Delta^1\searrow(\partial\Delta^m\boxtimes\Delta^1)\cup(\Delta^m\boxtimes 2)$.
By repeating the same process for smaller simplicies, we see that
$\Delta^m\boxtimes\Delta^1\searrow \Delta^m\boxtimes2$.
This can be applied to any simplicial complex $K$, and we see that $K\boxtimes \Delta^1$ is collapsible onto
$K\boxtimes2$.

\subsection{Massey product in $\Tor^{\fk[v_1,\cdots,v_m]}(\fk[K],\fk)$ }

Finally, we review the Massey products in $\Tor^{\fk[v_1,\cdots,v_m]}(\fk[K],\fk)$, according to
Section 3.2 of Buchstaber and Panov's book \cite{BP}. Here, we remark that we adopt 
the different conventions for the grading of $\Tor^{\fk[v_1,\cdots,v_m]}(\fk[K],\fk)$.

Recall that the torsion algebra $\Tor^{\fk[v_1,\cdots,v_m]}(\fk[K],\fk)$ is defined as the homology group of
the Koszul differential algebra $(\Lambda[u_1,\cdots,u_m]\otimes\fk[K],d)$,
where their bigrading and the differentials are defined by
$\deg u_i=(1,2)$, $\deg v_i=(0,2)$, and $du_i=v_i$, for $i=1,\cdots,m$.

A {\it multigrading} of the torsion algebra $\Tor^{\fk[v_1,\cdots,v_m]}(\fk[K],\fk)$ is defined by setting
\[
\mathrm{mdeg}\ u_1^{\varepsilon_1}\cdots u_m^{\varepsilon_m}v_1^{i_1}\cdots v_m^{i_m}=
(\varepsilon_1+\cdots+\varepsilon_m,2(i_1+\varepsilon_1),\cdots,2(i_m+\varepsilon_m)),
\]
and it is easy to see that
\[
\Tor^{\fk[v_1,\cdots,v_m]}_*(\fk[K],\fk)=\bigoplus_{\mathbf{a}\in\mathbb{N}^m}
\Tor^{\fk[v_1,\cdots,v_m]}_{*,2\mathbf{a}}(\fk[K],\fk).
\]

A subset $I\subset[m]$ may be viewed as a $(0,1)$-vector in $\mathbb{N}^m$, whose $i$-th coordinate
is $1$ if $i\in I$ and $0$ otherwise. Then, the following multigraded version of Hochster's formula holds.

\bigskip

{\bf Theorem 5.3.1(Theorem 3.2.9 of \cite{BP}).} {\it For any subset $I\subset[m]$, we have that
\[
\Tor^{\fk[v_1,\cdots,v_m]}_{i,2I}(\fk[K],\fk)\cong \tilde{H}^{|I|-i-1}(K_I;\fk),
\]
and $\Tor^{\fk[v_1,\cdots,v_m]}_{i,2\mathbf{a}}(\fk[K],\fk)=0$ unless $\mathbf{a}$ is a $(0,1)$-vector.
}
\bigskip

We recall the definition of the Massey product. Let $n$ be an integer greater than 1. 
Define $\bar{a}$ to be $(-1)^{\deg{a}+1}a$. Then, the $n$-fold Massey product 
$\langle a_1,\cdots,a_n\rangle$ is defined
to be the set of all homology classes represented by elements of the form
\[
a_{1,n}=\bar{a}_{1,1}a_{2,n}+\bar{a}_{1,2}a_{3,n}+\cdots+\bar{a}_{1,n-1}a_{n,n},
\]
for all solutions of the equations
\begin{eqnarray*}
a_i&=&[a_{i,i}]\quad \text{for $1\leq i\leq n$},\\
da_{i,j}&=&\bar{a}_{i,i}a_{i+1,j}+\bar{a}_{i,i+1}a_{i+2,j}+\cdots+\bar{a}_{i,j-1}a_{j,j}\quad
\text{for $1\leq i< j\leq n, (i,j)\neq (1,n)$}.\\
\end{eqnarray*}
By inductive argument it is easy to see that 
if $a_i\in \Tor^{\fk[v_1,\cdots,v_m]}_{k_i,2\mathbf{a}_i}(\fk[K],\fk)$, then \\
\[
a_{i,j}\in \Tor^{\fk[v_1,\cdots,v_m]}_{k_i+\cdots+k_j+(j-i),2(\mathbf{a}_i+\cdots+\mathbf{a}_j)}(\fk[K],\fk)\quad 
\text{for $1\leq i\leq j\leq n$, $(i,j)\neq(1,n)$},
\] 
and therefore, 
\[
\langle a_1,\cdots,a_n\rangle\ \subset 
\Tor^{\fk[v_1,\cdots,v_m]}_{k_1+\cdots+k_n+(n-2),2(\mathbf{a}_1+\cdots+\mathbf{a}_n)}(\fk[K],\fk).
\] 
The Massey product $\langle a_1,\cdots,a_n\rangle$ is said to be trivial if it contains the zero-element. 

\bigskip

{\bf Proposition 5.3.2.} {\it If an $n$-fold Massey product  $\langle a_1,\cdots,a_n\rangle$ is defined for \\
$a_i\in \Tor_{*,2I_i}^{\fk[v_1,\cdots,v_m]}(\fk[K],\fk)$, 
and $I_k\cap I_\ell\neq\emptyset$ for some $k\neq \ell$,
then $\langle a_1,\cdots,a_n\rangle$ is trivial. }

\bigskip

{\bf Proof.} Let $n\geq2$, and consider
an $n$-fold Massey product. As is observed above we have 
$\langle a_1,\cdots,a_n\rangle\ \subset \ \Tor_{*,2(I_1+\cdots+I_n)}^{\fk[v_1,\cdots,v_m]}(\fk[K],\fk).$
If $I_k\cap I_\ell\neq\emptyset$ for some $k\neq \ell$, then $I_1+\cdots+I_n$ is not a $(0,1)$-vector. Therefore, it follows from Theorem 5.3.1
that $\Tor_{*,2(I_1+\cdots+I_n)}^{\fk[v_1,\cdots,v_m]}(\fk[K],\fk)=0$. Thus, we complete the proof.
 \hfill $\Box$

\bigskip

{\bf Corollary 5.3.3.} {\it Let $K=K_1\cup_{\alpha}K_2$ be a simplicial complex on $[m]$, which is
obtained from two simplicial complexes $K_1$ and $K_2$ by gluing along a common simplex
$\alpha$ that is neither equal to $V(K_1)$ nor $V(K_2)$. Then, $K^*$ is non-Golod over a field $\fk$
if and only if there are faces $\sigma$ and $\tau$ of $K$ satisfying the following conditions:
\begin{enumerate}
\item[(1)] $V(K_1)-\alpha\subset\sigma$, $V(K_2)-\alpha\subset\tau$, and $\sigma\cup\tau=[m]$,
\item[(2)] the inclusion map 
\[
\st_{\link_K(\sigma\cap\tau)}(\sigma-\sigma\cap\tau)
\cup \st_{\link_K(\sigma\cap\tau)}(\tau-\sigma\cap\tau)\to \link_K(\sigma\cap\tau)
\]
induces a non-trivial map in the homology theory with coefficients $\fk$.
\end{enumerate}
}

\bigskip

{\bf Proof.} First, we show that Massey products for $n>2$ in  $\Tor^{\fk[v_1,\cdots,v_m]}_*(\fk[K^*],\fk)$
are trivial for any simplicial complex $K=K_1\cup_{\alpha}K_2$. This reduces the problem to proving that a 
non-trivial cup product exists if and only if the conditions (1) and (2) hold. 

 Let $a_i\in \Tor_{*,2I_i}^{\fk[v_1,\cdots,v_m]}(\fk[K^*],\fk)$ for $i=1,\cdots,n$, where $n>2$.
Then, we want to prove that the $n$-fold Massey product $\langle a_1,\cdots,a_n\rangle$ is trivial. If $J$ is a face of $K^*$, then \\
$\Tor_{*,2J}^{\fk[v_1,\cdots,v_m]}(\fk[K^*],\fk)\cong\tilde{H}^*((K^*)_J;\fk)=0$.
By this fact and the proposition above, we may assume that
$I_i=\sigma_i^c$ for some simplex $\sigma_i\in K$ for $i=1,\cdots,n$, and $I_i\cap I_j=\emptyset$ for $i\neq j$.
This implies that $\sigma_i\cup\sigma_j=[m]$. Thus, we may also assume that
$v_0\in I_1$, where $v_0\in V(K_1)-\alpha$.
Then, $v_0\not\in I_j=\sigma_j^c$, i.e., $v_0\in\sigma_j$ for $j>1$. Because a face of $K$
containing the vertex $v_0$ is a subset of $V(K_1)$, we see that
$V(K_2)-\alpha=V(K_1)^c\subset\sigma_j^c=I_j$ for $j>1$,
which contradicts the assumption that $I_j\cap I_k=\emptyset$ for $j\neq k>1$.
Thus, we have proved that Massey products for $n>2$ in 
$\Tor_*^{\fk[v_1,\cdots,v_m]}(\fk[K^*],\fk)$ are trivial.

It follows from Theorem 1.1 that $K^*$ is non-Golod if and only if there are disjoint subsets $I$ and $J$ of $[m]$ such that
\[
\iota_{I,J}:(K^*)_{I\sqcup J}=(\link_K(I^c\cap J^c))^*_{I\sqcup J}\to (K^*)_I*(K^*)_J=(\link_K(I^c))^*_I*(\link_K(J^c))^*_J
\] induces a non-trivial map in the cohomology theory. Then, $I^c$ and $J^c$ must be faces of $K$, and $I^c\cup J^c=[m]$, because $I\cap J=\emptyset$. Thus, we may assume that 
$V(K_1)-\alpha\subset I^c=\sigma$ and $V(K_2)-\alpha\subset J^c=\tau$. 
By the Alexander duality and its naturality, it follows that $\iota_{I,J}$
induces a non-trivial map in the cohomology theory if and only if the dual map
\begin{eqnarray*}
\iota_{I,J}^*:(\link_K(\sigma)^**\link_K(\tau)^*)^*_{I\sqcup J}
&=&\link_K(\sigma)*\Delta^J\cup \Delta^I*\link_K(\tau)\\
&=&\st_{\link_K(\sigma\cap\tau)}(\sigma-\sigma\cap\tau)
\cup \st_{\link_K(\sigma\cap\tau)}(\tau-\sigma\cap\tau)\\
&\to& \link_K(\sigma\cap \tau)
\end{eqnarray*}
induces a non-trivial map in the homology theory, where the first equality follows from Lemma 3.1.
Thus, we complete the proof.  \hfill $\Box$

\bigskip

{\bf Example 5.3.4.} Let $L$ be a simplicial complex on $\{2,3,4,5,6,7\}$ with facets
\[
267,\ 367,\ 467,\ 567,\ 236, \ 456, \ 2345.
\]
Put $\alpha=2345$ and consider the simplicial complex $K$ obtained from $\Delta^{\alpha+1}$ and $L$ 
by gluing along a common simplex $\alpha$. 
 Then, the Alexander dual of $K=\Delta^{\alpha+1}\cup_\alpha L$ is not Golod. In fact, for $\sigma=12345$ and
$\tau=67$, the inclusion map
\[
\st_K(\sigma)\cup\st_K(\tau)=\Delta^{\alpha+1}\cup \{\emptyset, 2,\ 3,\ 4,\ 5\}*\Delta^{\{6,7\}}\to K
\]
induces a non-trivial map in one dimensional homology groups. Needless to say, $\alpha$ is a facet of $L$.

\section{Proof of Theorem 2.7}

In this section, we prove Theorem 2.7. The proof is divided into three parts. In the first part, we show
that $K^*$ is Golod. Next, we show that $K^*$ is not stably homotopy Golod.
Finally, we demonstrate that $Z_{K^*}$ is torsion free for a particular $K$.

Before beginning the proof, we study the topology of the space $|(S^3_k\boxtimes \Delta^1)\cup_{\eta_k} S^2_4|$. As remarked at the end of Example 5.2.2,
there exists a sequence of elementary collapses
that collapses $S^3_k\boxtimes \Delta^1$ onto $S^3_k\boxtimes2=S^3_k$. 
This collapsion induces a deformation retraction $|S^3_k\boxtimes \Delta^1|$ onto $|S^3_k\boxtimes2|=|S^3_k|$, 
and the composite 
\[
|S^3_k|=|S^3_k\boxtimes1|\hookrightarrow |S^3_k\boxtimes \Delta^1|\to |S^3_k\boxtimes2|=|S^3_k|
\]
is homotopic to the identity map. This deformation retraction induces a deformation retraction 
$$\pi:|(S^3_k\boxtimes \Delta^1)\cup_{\eta_k} S^2_4|\to |S^2_4|$$ such that
$$\pi\circ |i|\simeq |\eta_k|:|S^3_k|\xrightarrow{|i|} |(S^3_k\boxtimes \Delta^1)\cup_{\eta_k} S^2_4|
\xrightarrow{\pi} |S^2_4|,$$ where $i$ is the inclusion map.

\subsection{$K^*$ is Golod}
Recall that
\[
K=\Delta^{V+v_0}\cup ((S^3_k\boxtimes \Delta^1\cup_{\eta_k}S^2_4)\cup F(S^3_k))*\Delta^{\{w_1\}}
\cup S^3_k*\Delta^{\{w_1,w_2\}},
\]
where $F(S^3_k)$ is the simplicial complex defined in \S2. We apply Corollary 5.3.3 to  
$K=K_1\cup_\alpha K_2$, where  
\[
K_1=\Delta^{V+v_0}, \quad K_2=\Delta^V\cup ((S^3_k\boxtimes \Delta^1\cup_{\eta_k}S^2_4)\cup F(S^3_k))*\Delta^{\{w_1\}}\cup S^3_k*\Delta^{\{w_1,w_2\}}, 
\]
$\alpha=V$, $V(K_1)=V+v_0$, and $V(K_2)=V\sqcup \{w_1,w_2\}$.  
Thus, by Corollary 5.3.3 we only have to prove that the map
\[
\st_{\link_K(\sigma\cap\tau)}(\sigma-\sigma\cap\tau)
\cup \st_{\link_K(\sigma\cap\tau)}(\tau-\sigma\cap\tau)\to \link_K(\sigma\cap\tau)
\]
induces the trivial map in the homology theory,
where $\sigma$ and $\tau$ are faces of $K$ such that $v_0\in\sigma$, $\{w_1,w_2\}\subset\tau$,
and $\sigma\cup\tau=V(K)$.

First, we consider the case with $\sigma\cap\tau=\emptyset$. That is, we prove that the map
\[
\st_K(\sigma)\cup \st_K(\tau)=\Delta^{V+v_0}\cup \st_K(\tau)\to K
\]
induces the trivial map in the homology theory. Because this map factors as
\[
\Delta^{V+v_0}\cup \st_K(\tau)\to\Delta^{V+v_0}\cup S^3_k*\Delta^{\{w_1,w_2\}}\to K,
\]
it suffices to prove that the inclusion map
\[
\Delta^{V+v_0}\cup S^3_k*\Delta^{\{w_1,w_2\}}\to K
\]
induces the trivial map in the homology theory.
Because $|\Delta^{V+v_0}\cup S^3_k*\Delta^{\{w_1,w_2\}}|\simeq\Sigma|S^3_k|=S^4$ and
$|K|\simeq \Sigma|S^2_4|=S^3$, the above inclusion map induces the trivial map in the homology theory.

Next, we consider the case that $\rho=\sigma\cap\tau$ is a non-empty face of $S^3_k$. 
Let $L=L_1\cup_\alpha L_2$ be a simplicial complex obtained from $L_1$ and $L_2$ by gluing along a 
common face $\alpha$, and let $\beta$ be a face of $\alpha$. Then $\link_L(\beta)=\link_{L_1}(\beta)
\cup_{\alpha-\beta}\link_{L_2}(\beta)$. By this observation we see that  
\[
\link_K(\rho)
=\Delta^{(V+v_0)-\rho}\cup (\link_{S^3_k\boxtimes \Delta^1\cup_{\eta_k} S^2_4}(\rho)
\cup F(\link_{S^3_k}(\rho))*\Delta^{w_1}\cup\link_{S^3_k}(\rho)*\Delta^{\{w_1,w_2\}}. 
\]
The situation is, therefore, similar to the case with $\rho=\emptyset$. Here, we are only required 
to prove that the map
\[
f:|\link_{S^3_k}(\rho)|\to |\link_{S^3_k\boxtimes \Delta^1\cup_{\eta_k} S^2_4}(\rho)|
\]
induced by the inclusion induces the trivial map in the homology theory. $f$ factors through a contractible space
$|\link_{S^3_k\boxtimes \Delta^1}(\rho)|$.
In fact, because $S^3_k\boxtimes \Delta^1$ is a triangulation
of $S^3\times I$ and $\rho$ is a non-empty simplex in the boundary,
it follows that $\link_{S^3_k\boxtimes \Delta^1}(\rho)$ is a triangulation of a hemisphere.
Thus, $f$ clearly induces the trivial map in the homology theory, and we complete the proof.  \hfill $\Box$

\subsection{$K^*$ is not stably homotopy Golod}

Next, we show that the map
\begin{equation}
|\iota_{\{w_1,w_2\},V+v_0}|:
|K^*|\to |(K^*)_{\{w_1,w_2\}}*(K^*)_{V+v_0}|
\end{equation}
is stably non-trivial. To show this, we consider the mapping cone of $|\iota_{\{w_1,w_2\},V+v_0}|$,
which is denoted by $C_{|\iota_{\{w_1,w_2\},V+v_0}|}$, and show that $Sq^2$ acts
non-trivially on its mod-2 cohomology groups. Let $m=|V+\{v_0,w_1,w_2\}|$. Then 
we have the following isomorphisms of (co)homology groups:
\begin{eqnarray*}
\tilde{H}^p(C_{|\iota_{\{w_1,w_2\},V+v_0}|})&\cong&
H^p( |(K^*)_{\{w_1,w_2\}}*(K^*)_{V+v_0}|,|K^*|)\\
&\cong&H_{m-2-p}(K, \Delta^{V+v_0}\cup  S^3_k*\Delta^{\{w_1,w_2\}})\\
&\cong&H_{m-2-p}(\Delta^{V+v_0}\cup(S^3_k\boxtimes\Delta^1\cup_{\eta_k}S^2_4)*\Delta^{\{w_1\}},
\Delta^{V+v_0}\cup S^3_k*\Delta^{\{w_1\}})\\
&\cong&H_{m-2-p}(\Sigma|S^3_k\boxtimes\Delta^1\cup_{\eta_k}S^2_4|, \Sigma|S^3_k|)\\
&\cong&\tilde{H}_{m-2-p}(\Sigma\mathbb{C}P^2),
\end{eqnarray*}
where the second isomorphism follows from the Alexander duality and the fact that 
\begin{eqnarray*}
((K^*)_{\{w_1,w_2\}}*(K^*)_{V+v_0})^*&=&  ((K^*)_{\{w_1,w_2\}})^*_{\{w_1,w_2\}}*\Delta^{V+v_0}\cup\Delta^{\{w_1,w_2\}}*((K^*)_{V+v_0})^*_{V+v_0})\\
&=&\link_K(V+v_0)*\Delta^{V+v_0}\cup\Delta^{\{w_1,w_2\}}*\link_K(\{w_1,w_2\})\\
&=&\Delta^{V+v_0}\cup S^3_k*\Delta^{\{w_1,w_2\}}.
\end{eqnarray*}
In the third isomorphism we deformed $S^3_k*\Delta^{\{w_1,w_2\}}$ onto $S^3_k*\Delta^{\{w_1\}}$, and 
after that we deformed $F(S^3_k)*\Delta^{\{w_1\}}$ onto $S^3_k*\Delta^{\{w_1\}}$, in the fourth 
isomorphism we collapses $|\Delta^{V+v_0}|$ to the point, and in the fifth isomorphism we used 
the fact that $|(S^3_k\boxtimes \Delta^1)\cup_{\eta_k} S^2_4|$ is the mapping cylinder of 
$\eta:|S^3_k|\to |S^2_4|$ and $\mathbb{C}P^2\simeq S^2\cup_\eta e^4$ where $e^4$ denotes a 
4-dimensional cell.  
Since $Sq^2_*$ acts non-trivially on $\tilde{H}_5(\Sigma\mathbb{C}P^2)$,  
it follows from Lemma 5.1.1 that $Sq^2$ acts non-trivially on 
$\tilde{H}^{m-7}(C_{|\iota_{\{w_1,w_2\},V+v_0}|})$, which implies that the map (6.1) is stably non-trivial. \hfill $\Box$

\subsection{$Z_{K^*}$ is torsion free}

Finally, we will demonstrate that $Z_{K^*}$ is torsion free if $K$ is constructed using the simplicial map
$\eta_{12}:S^3_{12}\to S^2_4$ described in \cite{MS}. We believe that $Z_{K^*}$ is torsion free
in general, but we are currently unable to prove this.

By Theorem 1.1 and the Alexander duality, we have the following isomorphisms:
\[
H^p(Z_{K^*};\Z)\cong \bigoplus_{I\subset[m]}\tilde{H}^{p-|I|-1}((K^*)_I:\Z)\cong\bigoplus_{I\subset[m]}\tilde{H}_{2|I|-p-2}(\link_K(I^c):\Z),
\]
where we have used the fact that $(K^*)_I=(\link_K(I^c))^*$. Thus, to show that $Z_{K^*}$ is torsion free, it suffices to show that
every $\link_K(\sigma)$ is torsion free for all faces $\sigma$ of $K$.
The longest part of the computation is the following, and the remaining parts are omitted.

\bigskip

{\bf Lemma 6.3.1.} {\it $\link_{S^3_{12}\boxtimes \Delta^1\cup_{\eta_{12}} S^2_4}(\sigma)$ is torsion free
for any simplex $\sigma$ in $S^3_{12}\boxtimes \Delta^1\cup_{\eta_{12}} S^2_4$.}

\bigskip

{\bf Proof.} Set $L=S^3_{12}\boxtimes \Delta^1\cup_{\eta_{12}} S^2_4$.
If $\sigma=\emptyset$, then clearly $\link_L(\sigma)=L$ is torsion free.

Let $v$ be a vertex of $S^3_{12}$. The Mayer-Vietoris sequence associated with the decomposition
$L=(L-v)\cup \st_L(v)$ reduces to the long exact sequence
\[
\cdots\to \tilde{H}_i(\link_L(v))\to \tilde{H}_i(L-v)\to \tilde{H}_i(L)\to \cdots.
\]
The sequence of elementary collapses of
$S^3_{12}\boxtimes \Delta^1$ onto $S^3_{12}\boxtimes2$ given in Example 5.2.2 induces a sequence of
elementary collapses of $S^3_{12}\boxtimes\Delta^1\cup_{\eta_{12}} S^2_4$ onto $S^2_4$. Moreover,
this sequence of elementary collapses induces a collapse of $L-v$ onto $S^2_4$. This means that
the inclusion map $|L-v|\to |L|$ is homotopy equivalent, and that $\tilde{H}_*(\link_L(v))=0$.

To proceed further with the computation, we need to know the concrete structure of $L$. $S^3_{12}$ has
twelve vertices $\{a_i, b_i,  c_i, d_i\}_{i=0,1,2}$, and the following is the list of its facets:

\[
\begin{matrix}
a_0b_0c_0c_1&a_0b_0b_1c_1&a_0a_1b_1c_1&a_1a_2b_1c_1&a_2b_1c_1c_2&a_2b_1b_2c_2&
a_2b_0b_2c_2&a_0a_2b_0c_2&a_0b_0c_0c_2\\
a_0a_2b_0d_1&a_0b_0b_1d_1&b_0b_1c_1d_1&
b_1c_1c_2d_1&a_2c_1c_2d_1&a_0a_2c_2d_1&a_0b_1d_0d_1&b_1c_2d_0d_1&a_0c_2d_0d_1\\
a_0b_1d_0d_2&b_1c_2d_0d_2&a_0c_2d_0d_2&a_0a_1b_1d_2&a_1a_2b_1d_2&a_2b_1b_2d_2&
a_2b_0b_2d_2&b_1b_2c_2d_2&b_0b_2c_2d_2\\
b_0c_0c_2d_2&b_0c_0c_1d_2&a_0c_0c_2d_2&
a_0c_0c_1d_2&a_0a_1c_1d_2&a_1a_2c_1d_2&a_2b_0d_1d_2&b_0c_1d_1d_2&a_2c_1d_1d_2\\
\end{matrix}
\]

We order the vertices as $a_0<a_1<a_2<b_0<\cdots<c_0<\cdots<d_0<d_1<d_2$. The vertices of
$S^2_4$ are $a,b,c,d$, and are ordered such that $a<b<c<d$.
Furthermore, $\eta_{12}:S^3_{12}\to S^2_4$ is the map defined by the mappings $a_i\mapsto a$,
$b_i\mapsto b$, $c_i\mapsto c$, and $d_i\mapsto d$, for $i=0,1,2$. Then, the list of the
facets of $L$ is as follows:
\[
\begin{matrix}
a_0b_0c_0c_1c&&&a_0a_2b_0d_1d&a_0a_2b_0bd&a_0a_2abd&a_0b_1d_0d_2d&&\\
a_0b_0b_1c_1c&a_0b_0b_1bc&&a_0b_0b_1d_1d&a_0b_0b_1bd&&b_1c_2d_0d_2d&&\\
a_0a_1b_1c_1c&a_0a_1b_1bc&a_0a_1abc&b_0b_1c_1d_1d&b_0b_1c_1cd&b_0b_1bcd&a_0c_2d_0d_2d&&\\
a_1a_2b_1c_1c&a_1a_2b_1bc&a_1a_2abc&b_1c_1c_2d_1d&b_1c_1c_2cd&&a_0a_1b_1d_2d&a_0a_1b_1bd&a_0a_1abd\\
a_2b_1c_1c_2c&&&a_2c_1c_2d_1d&a_2c_1c_2cd&&a_1a_2b_1d_2d&a_1a_2b_1bd&a_1a_2abd\\
a_2b_1b_2c_2c&a_2b_1b_2bc&&a_0a_2c_2d_1d&a_0a_2c_2cd&a_0a_2acd&a_2b_1b_2d_2d&a_2b_1b_2bd&\\
a_2b_0b_2c_2c&a_2b_0b_2bc&&a_0b_1d_0d_1d&&&a_2b_0b_2d_2d&a_2b_0b_2bd&\\
a_0a_2b_0c_2c&a_0a_2b_0bc&a_0a_2abc&b_1c_2d_0d_1d&&&b_1b_2c_2d_2d&b_1b_2c_2cd&b_1b_2bcd\\
a_0b_0c_0c_2c&&&a_0c_2d_0d_1d&&&b_0b_2c_2d_2d&b_0b_2c_2cd&b_0b_2bcd\\
&&&&&&b_0c_0c_2d_2d&b_0c_0c_2cd&\\
&&&&&&b_0c_0c_1d_2d&b_0c_0c_1cd&\\
&&&&&&a_0c_0c_2d_2d&a_0c_0c_2cd&\\
&&&&&&a_0c_0c_1d_2d&a_0c_0c_1cd&\\
&&&&&&a_0a_1c_1d_2d&a_0a_1c_1cd&a_0a_1acd\\
&&&&&&a_1a_2c_1d_2d&a_1a_2c_1cd&a_1a_2acd\\
&&&&&&a_2b_0d_1d_2d&&\\
&&&&&&b_0c_1d_1d_2d&&\\
&&&&&&a_2c_1d_1d_2d&&\\
\end{matrix}
\]
We see that $L=L_1\cup L_2$, where $L_1$ is the subcomplex generated by the
facets on the first to third columns from the left in the above list, and $L_2$ is the subcomplex
generated by the other facets. Then, $L=\st_L(c)\cup \st_L(d)$.
If $\sigma\cap\{c,d\}=\emptyset$, then $\link_L(\sigma)=\link_{L_1}(\sigma+c)*c\cup\link_{L_2}(\sigma+d)*d$,
which implies that $|\link_L(\sigma)|$ is homotopy equivalent to a suspension space. Therefore,
if $\link_L(\sigma)\leq2$, i.e., $\dim\sigma\geq1$, then it follows that $\link_L(\sigma)$ is torsion free, by
dimensional reasoning. We consider a face with $\dim\sigma\leq0$. For a vertex or the empty face in $S^3_{12}$,
we have already proved that $\link_L(\sigma)$ is torsion free. If $v=a$ or $v=b$, then we see that
$|\link_L(v)|$ is homotopy equivalent to $S^3$. Thus, we have completed the proof in this case.

If $c\in\sigma$ or $d\in\sigma$, then we see by direct computation that
$|\link_L(c)|$ and $|\link_L(d)|$ are homotopy equivalent to $S^3$. Moreover,
$\link_L(ac)$, $\link_L(bc)$, $\link_L(cd)$, $\link_L(xc)$, $\link_L(ad)$,
$\link_L(bd)$, and $\link_L(xd)$ are triangulations of $S^2$, where $x$ is a vertex of $S^3_{12}$.
If $\dim\link_L(\sigma)\leq1$, then $\link_L(\sigma)$ is torsion free by
dimensional reasoning, and the proof is complete for all cases. \hfill $\Box$

\end{document}